\documentclass{AIMS}
\usepackage{amsmath}
  \usepackage{paralist}
  \usepackage{graphics} 
  \usepackage{epsfig} 
\usepackage{graphicx}  \usepackage{epstopdf}
 \usepackage[colorlinks=true]{hyperref}
\hypersetup{urlcolor=blue, citecolor=red}

\def\currentvolume{X}
 \def\currentissue{X}
  \def\currentyear{200X}
   \def\currentmonth{XX}
    \def\ppages{X--XX}
     \def\DOI{10.3934/xx.xx.xx.xx}

\newcommand{\bR}{{\mathbb{R}}}

\newtheorem{theorem}{Theorem}[section]

\newtheorem{proposition}{Proposition}

\theoremstyle{definition}
\newtheorem{definition}[theorem]{Definition}

\title[EBT for a system of age-structured equations] 
      {The Escalator Boxcar Train method for a system of age-structured equations}

\author[Piotr Gwiazda, Karolina Kropielnicka, Anna Marciniak-Czochra]{}


 \email{pgwiazda@mimuw.edu.pl}
 \email{karolina.kropielnicka@mat.ug.edu.pl}
\email{anna.marciniak@iwr.uni-heidelberg.de}

\begin{document}
\maketitle

\centerline{\scshape Piotr Gwiazda }
\medskip
{\footnotesize
 \centerline{Institute of Applied Mathematics and Mechanics, }
   \centerline{University of Warsaw,}
   \centerline{ Warszawa 02-097, Poland}
} 

\medskip

\centerline{\scshape Karolina Kropielnicka}
\medskip
{\footnotesize
 \centerline{ Institute of Mathematics}
   \centerline{University of Gda\'nsk}
   \centerline{80-952 Gda\'nsk, Poland}
}

\medskip

\centerline{\scshape Anna Marciniak-Czochra}
\medskip
{\footnotesize
 \centerline{ Institute of Applied Mathematics, Interdisciplinary Center for scientific Computing (IWR) and BIOQUANT}
   \centerline{University of Heidelberg}
   \centerline{69120 Heidelberg, Germany}
}

\bigskip

 \centerline{(Communicated by the associate editor name)}

\begin{abstract}
The Escalator Boxcar Train method (EBT) is a numerical method for structured population models of McKendrick -- von Foerster type. Those models consist of a certain class of hyperbolic partial differential equations and describe time evolution of the distribution density of the structure variable describing a feature of individuals in the population. The method was introduced in late eighties and widely used in theoretical biology, but its convergence was proven only in recent years using the framework of measure-valued solutions. Till now the EBT method was developed only for scalar equation models. In this paper we derive a full numerical EBT scheme for age-structured, two-sex population model (Fredrickson-Hoppensteadt model), which consists of  three coupled hyperbolic partial differential equations with nonlocal boundary conditions. It is the first step towards extending the EBT method to systems of structured population equations. 

\end{abstract}

\section{Introduction}
\subsection{Problem formulation}

The Escalator Boxcar Train (EBT) algorithm is a numerical method for solving structured population models. It is based on representing the solution as a sum of masses localised in discrete points and tracing their evolution due to transport and growth. The algorithm has been used in applications for a long time \cite{deRoos89}, however convergence of the scheme was shown only recently in ref.~\cite{EBT} (without the rate of the convergence), and later in ref.~\cite{GwJabMarUli14} (the rate of the convergence) using the metric space  approach proposed in ref.~\cite{GLMC,GMC}. So far the method has been established only for single equation models. The aim of this paper is to derive the EBT algorithm for a system of structured population equations.

We focus on the Fredrickson-Hoppensteadt model, which is a two-sex population model with age as a structure variable. The model was originally formulated in ref.~\cite{Fredrickson71} and later developed in ref.~\cite{Hoppensteadt75}. The model consists of three population equations with structure coupled through nonlocal boundary terms and a nonlocal and nonlinear source term. Functions $u^m(t,x)$ and  $u^f(t,x)$ describe the distribution of males and females, while $u^c(t,x,y)$ is the number of couples at time $t$. Structural variables $x$ and $y$ denote the age of males and females, respectively. The following system of equations describes the dynamics of the population of males, females and couples:

\begin{eqnarray}
\partial_t u^m(t,x)+\partial_x u^m(t,x)+c^m(t,x)u^m(t,x)&=&0,\nonumber\\
u^m(t,0)&=&\int_{\bR_+^2}\beta^m(t,x,y)u^c(t,x,y)dxdy,\nonumber\\
u^m(0,x)&=& u^m_0(x),\nonumber\\\nonumber
\\\nonumber
\partial_t u^f(t,y)+\partial_y u^f(t,y)+c^f(t,y)u^f(t,y)&=&0,\nonumber\\
u^f(t,0)&=&\int_{\bR_+^2}\beta^f(t,x,y)u^c(t,x,y)dxdy,\nonumber\\
u^f(0,x)&=& u^f_0(x),\nonumber\\ \label{the_system}
\end{eqnarray}
\begin{eqnarray*}
\partial_t u^c(t,x,y)+\partial_{x} u^c(t,x,y)+\partial_{y}u^c(t,x,y)+c^c(t,x,y)u^c(t,x,y)&=&\mathcal{T}(t,x,y),\\
u^c(t,x,0)=u^c(t,0,y)&=&0,\\
u^c(0,x,y)&=& u^c_0(x,y).\\
\end{eqnarray*}

Functions $c^m$, $c^f$ and $c^c$ describe the rate of disappearance of individuals. In the case of males and females disappearance equals death, while in the case of couples it equals divorce or death of one of spouses. Functions $\beta^m$ and $\beta^f$ are birth rates of males and females. In the general case the functions may depend nonlinearly on ecological pressure, which is usually modelled by a nonlocal operator depending on the distribution of males, females and couples. However, in this paper we restrict our considerations to the case where disappearance and birth rates depend only on time and on the structure variable. 

The marriage function $\mathcal{T}$ models the number of new marriages of males and females of age $x$ and $y$, respectively, at time $t$. It  depends nonlinearly on the distribution of individuals. The choice of this function is a subject of ongoing discussions, see \cite{Hadeler89}, \cite{Hadeler88},  \cite{Pruss94}, \cite{Martcheva99},  due to the properties like heterosexuality, homogeneity, consistency or competition. In this paper we follow the formulation proposed in ref.~\cite{Inaba93}, namely:

\begin{equation}\label{marriage_function}
\mathcal{T}(t,x,y)=F(t,u^m(t,x),u^f(t,x),u^c(t,x,y))=
\end{equation}
$$
\frac{\Theta(x,y)h(x)g(y)\left[u^m(t,x)-\int_0^\infty u^c(t,x,y)dy\right]\left[u^f(t,y)-\int_0^\infty u^c(t,x,y)dx\right]}{\gamma+\int_0^\infty h(x)\left[u^m(t,x)-\int_0^\infty u^c(t,x,y)dy\right]dx+\int_0^\infty g(y)\left[u^f(t,y)-\int_0^\infty u^c(t,x,y)dx\right]dy}.
$$

Function $\Theta(x,y)$ describes the marriage rate of males of age $x$ and females of age $y$. Notice that $\left[u^m(t,x)-\int_0^\infty u^c(t,x,y)dy\right]$ is the amount of unmarried males and $\left[u^f(t,y)-\int_0^\infty u^c(t,x,y)dx\right]$ is the number of unmarried females.  Function $h\in{\mathbf{L^1}}(\bR_+)\cap{\mathbf{L^\infty}}(\bR_+)$ describes the distribution of eligible males on the marriage market. Function $g$ is of the same regularity and describes the distribution of eligible females on the marriage market.

\subsection{EBT method for a scalar equation}\label{scalarproblemanalysis}
Before moving to the derivation of the numerical scheme for system (\ref{the_system}), we  briefly discuss the analytical and numerical results for a simpler model of the McKendrick - von Foerster type \cite{McKendrick26}:

\begin{eqnarray}\label{age_structured}
\partial_t u(t,x)+\partial_x (b(t,x)u(t,x))+c(t,x)u(t,x)&=&0,\\
u(t,x_b)&=&\int_{x_b}^\infty\beta(t,x)u(t,x)dx,\nonumber\\
u(0,x)&=& u_0(x),\nonumber\\\nonumber
\end{eqnarray}
which is a one-sex population model with structure variable $x$. The function $b(t,x)$ is the rate at which the individuals change their state. In particular case, when $x$ represents age, function $b(t,x)$ equals one.
One of the many numerical methods applied to this problem is  the Escalator Boxcar Train introduced in ref.~\cite{deRoos89}. 

EBT method is based on representing the solution as a sum of masses (cohorts) localised in discrete points and tracing their spatio--temporal evolution using the following algorithm:

\begin{eqnarray}\label{EBT_org_system}
\left\{
	\begin{array}{rcll}
\displaystyle
\mathrm{\frac{d}{dt}} x_i(t) &=& b(t, x_i(t)),&
\mathrm{for}\;\; i=B+1, \dots, J,
\\[2mm]
\displaystyle
\mathrm{\frac{d}{dt}} m_i(t) &=& -c(t,x_i(t))m_i(t),
&
\mathrm{for}\;\; i=B+1, \dots, J,
	\end{array}
\right.
\end{eqnarray}

\begin{eqnarray}\label{EBT_org_BC}
\left\{
	\begin{array}{rcl}
\displaystyle
x_B(t) &=& 
\left\{
\begin{array}{lcl}
\frac{\pi_B(t)}{m_B(t)} + x_b,
&&
\mathrm{if}\;\; m_B(t) > 0,
\\
x_b, && \mathrm{otherwise},
\end{array}
\right.\\ \\
\displaystyle
{\mathrm{\frac{d}{dt}} \pi_B(t)} &=& b(t, x_b)m_B(t)
+ \partial_x b( t,x_b)\pi_B(t)
\\
&& - c( t,x_b)\pi_B(t),
\\[3mm]
\displaystyle
\mathrm{\frac{d}{dt}} m_B(t) &=&-c(t, x_b)m_B(t)  - \partial_x c (t, x_b)\pi_B(t)
\\
&& + \sum_{i=B}^{J} \beta(t, x_i(t))m_i(t),\\
m^B (t_k) &=& 0,\\
\pi^B (t_k) &=& 0.
	\end{array}
\right.
\end{eqnarray} 
 
 The procedure consists in solving the system of ordinary differential equations (ODEs) (\ref{EBT_org_system})-(\ref{EBT_org_BC}) on  a sufficiently short time interval $[t_k,t_{k+1}]$, iteratively with respect to $k=0,1,\ldots$. The details of the EBT algorithm are discussed in Section 2 for system (\ref{the_system}) and therefore we omit them at this stage.

As mentioned before, the rate of convergence of the method for a nonlinear version of model \eqref{age_structured} has been recently shown in ref.~\cite {GwJabMarUli14}. The results are based  on Lipschitz semiflows in metric spaces developed in ref.~\cite{ColomboGuerra2009} and applied to structured population models in ref.~\cite{GMC}.
\begin{definition}\label{DefLipschSemiflow}
Let $(E,\rho)$ be a metric space. A Lipschitz semiflow is a semigroup $S:[0,\delta]\times[0,T]\times E\rightarrow E$ satisfying
$$
\rho(S(t;\tau)\mu,S(s;\tau)\nu)\leq L(\rho(\mu,\nu)+|t-s|),
$$
where $s,t\in[0,\delta],\tau,s+\tau,t+\tau\in[0,T]$ and $\mu,\nu\in E$.
\end{definition}

The idea of the proof of convergence presented in ref.~\cite{GwJabMarUli14} consists in considering the system (\ref{age_structured}) (under the assumption that $b,c,\beta \in \mathbf C^{\alpha}([0,T];\mathbf{W}^{1,\infty}(\bR_+)$) as a semiflow in a space of nonnegative Radon measures (${\mathcal M^+}({\bR_+})$) equipped with the Lipschitz bounded distance, defined below.

\begin{definition}\label{defflatmetric} 
Let $\mu, \nu \in{\mathcal M^+}({\bR_{+}})$.  The distance
function $\rho_F : {\mathcal M^+}({\bR_{+}}) \times  {\mathcal M^+}({\bR_{+}}) \rightarrow
[0, \infty]$ is defined by
\begin{equation}\label{Flatmetric}
  \rho_F(\mu, \nu) :=
 \sup \Big\{\int_{\bR_{+}} \psi  d (\mu - \nu)
  \big| \: \psi \in \mathbf{C}^1({\bR_{+}}), \|\psi\|_{W^{1,\infty}} \leq 1 \Big\},
  \end{equation}
where $\|\psi\|_{W^{1,\infty}} = \max\{\|\psi\|_{\infty},\|\partial_x \psi\|_{\infty}\}.$
\end{definition}
\noindent

The metric $\rho_F$ is a distance derived from the dual norm of
$W^{1,\infty}$, with the space $W^{1,\infty}$ equipped with its usual norm, i.e.
$$
\|\gamma\|_{W^{1,\infty}} = \max\{ \|\gamma\|_{L^{\infty}}, \|\partial_x \gamma\|_{L^{\infty}}\},
$$
see ref.~\cite{MullerOrtiz} and \cite{Zhidkov}.\\

 The solution to the system (\ref{EBT_org_system})-(\ref{EBT_org_BC}) is interpreted as an element of the space of nonnegative Radon measures given by the formula 
\begin{equation}\label{measure_sol}
\mu^n(t)=\sum_{B(t)}^Jm_i(t)\delta_{\{x_i(t)\}},
\end{equation}
which means that $\mu^n(t)$ is a linear combination of Dirac deltas with masses localised in points $x_i(t)$. 

The proof of convergence of the ETB scheme is based on the following proposition, the proof of which is a modification of the proof presented in ref.~\cite{bressan_00}.
\begin{proposition}\label{tangential_ineq}
Let $S:E\times[0,\delta]\times[0,T]\rightarrow E$ be a Lipschitz semiflow. For every Lipschitz continuous map $[0,T]\ni t\mapsto\nu_t\in(E,\rho)$ the following estimate holds:
\begin{equation}\label{estimate}
\rho(\nu_t,S(t;0)\mu_0)\leq L\int_{[0,T]}\liminf_{h\rightarrow 0}\frac{\rho(\nu_{t+h},S(h;\tau)\nu_t)}{h}d\tau,
\end{equation}
where $\rho$ is a corresponding metric.
\end{proposition}
Application of Proposition \ref{tangential_ineq} reduces the proof of the convergence of the measure given by formula (\ref{measure_sol}) to the solution of the model (\ref{age_structured}) to the estimate of the right hand side of inequality (\ref{estimate}). In ref.~\cite{GwJabMarUli14} it was proven analytically and confirmed numerically that the method is of the first order.


\subsection{Organisation of the paper}
This paper is devoted to derivation of a corresponding EBT method for the two-sex age-structured population model \eqref{the_system}. The convergence of the method follows from stability results obtained in ref.~\cite{Ulikowska2012} and its proof is deferred to a forthcoming paper.  Section 2 is devoted to derivation of the EBT algorithm for the two-sex population model. In Section 3, we provide well-posedness results for the structured population model and for the EBT scheme, and discuss a sketch of the convergence proof.

\section{Derivation of the EBT algorithm for the two-sex population model}

This section is devoted to the formal derivation of the EBT scheme for system (\ref{the_system}).  We assume that the solution $u^f,u^m$ i $u^c$ to system (\ref{the_system}) are continuously differentiable with respect to time and twice time continuously differentiable with respect to the structure variables. Additionally, the solution itself and its first derivative with respect to time and first and second derivatives with respect to the structure variables are assumed to be bounded. The same regularity is assumed on all model coefficients appearing in the underlying system (\ref{the_system}).
 
The first step of the EBT method is grouping the initial distribution of males, females  and couples into the so-called cohorts, with the grouping performed with respect to the age of the individuals. While in the case of males and females the groups form one-dimensional intervals including individuals of certain age, in the case of couples we deal with two-dimensional intervals to take into account the age of both partners. Each cohort  is characterised by two quantities: the mass and its location. The mass ($m^m(t),\ m^f(t),\  m^c(t)$) is the number of individuals or couples within the cohort at time $t$, and its location ($x^m(t),\ y^f(t)$ and  $(x^c(t),y^c(t))$, respectively) is the average value of the structure variable over the underlying cohort. Cohorts are evolving in time along the characteristic lines of the underlying problem. We wish to track how masses and locations of cohorts change in time.  Once individuals or couples are assigned to a certain cohort they remain there till disappearance (death or split-up in case of couples). The influx of new couples is described with marriage function, while the influx of new males and females is given in the boundary conditions.  Each time step can be treated as an internalisation moment, where a new boundary cohort appears. In case of couples the boundary cohorts are empty and in case of individuals the masses in boundary cohorts result from the boundary conditions.

Following this argument, we start with introducing the initial cohorts. At time $t=0$ we impose space mesh points $l_i^m(0)$ and $l_j^f(0)$, where $i,j=B_0,\ldots,J$ in such a way that $l_{B_0}^m(0)=l_{B_0}^f(0)=0$. We define $J-B_0$ initial intervals for males and females

$$
\Omega^m_{i+1}(0)=\left[l_i^m(0),l_{i+1}^m(0)\right),\ \Omega^f_{j+1}(0)=\left[l_j^f(0),l_{j+1}^f(0)\right),\  i,j=B_0,\ldots,J-1,
$$
and $(J-B_0)^2$ two-dimensional intervals for couples
$$
\Omega^c_{(i+1)(j+1)}(0)=\Omega^m_{i+1}(0)\times \Omega^f_{j+1}(0),\ i,j=B_0,\ldots,J-1.
$$
The above partitioning requires the condition that the supports of initial functions  $u_0^m$ , $u_0^f$ and $u_0^c$, which are the initial distributions of individuals and couples, are contained in the sums of all corresponding cohorts that is
$$
{\rm supp}(u_0^m)\subset \Omega^m(0)=\bigcup_{i=B_0+1}^J \Omega_i^m(0),\ {\rm supp}(u_0^f)\subset \Omega^f(0)=\bigcup_{j=B_0+1}^J \Omega_j^f(0),\ 
$$
$$
{\rm supp}(u_0^c)\subset \Omega^c(0)=\bigcup_{i,j=B_0+1}^J \Omega_{ij}^c(0).
$$

\noindent Bearing in mind that the cohorts evolve in time, we write
\{$\Omega_i(t)\}^J_{i=B_0+1}$, \newline $\{\Omega^f_j(t)\}^J_{j=B_0+1}$, $\{\Omega^c_{ij}(t)\}^J_{i,j=B_0+1}$ and understand that functions $l_i^m(t)$ and $l_i^f(t)$ are characteristics defined by the transport operators, thus they satisfy differential equations $\frac{d}{dt}l_i^m(t)=1$ and $\frac{d}{dt}l_j^f(t)=1$ and are straight lines:
\begin{eqnarray*}
l_i^m(t)&=&t+l_i^m(0),\ i=B_0,\ldots,J,\\
l_j^f(t)&=&t+l_j^f(0),\ j=B_0,\ldots,J.\\
\end{eqnarray*}
We impose a mesh on the time variable $t\in[0,T)$ in such a way that $t_0=0$ and $\bigcup_{n=0}^N[t_n,t_{n+1})=[0,T)$. For $0\leq t < t_1$, we define boundary cohorts $\Omega_{B_0}^m(t)=[0,l_{B_0}^m(t))$, $\Omega_{B_0}^f(t)=[0,l_{B_0}^f(t))$ and $\Omega_{B_0B_0}^c(t)=\Omega_{B_0}^m(t)\times \Omega_{B_0}^f(t)$, $\Omega_{iB_0}^c(t)=\Omega_i^m(t)\times \Omega_{B_0}^f(t)$, $\Omega_{B_0j}^c(t)=\Omega_{B_0}^m(t)\times \Omega_j^f(t)$, where $i,j=B_0+1,J$. Mesh points $t_n$, where $t_n=1,\ldots,N$, are called internalisation moments.

Altogether, we have $J-B_0+1$ initial cohorts for males and females and $(J-B_0+1)^2$ initial cohorts of couples.

At each point $t_n$, we define $B_n=B_0-n$ and $l_{B_n}(t)=t-t_n$ for $t_n\leq t<t_{n+1}$. Now, for $t \in [t_n,t_{n+1})$, we define internal cohorts: 
\begin{eqnarray*}
\Omega_{i+1}^m(t)&=&[l_i^m(t),l_{i+1}^m(t)),\ {\rm for}\ i=B_{n-1}‚Äô\ldots,J-1\\
\Omega_{j+1}^f(t)&=&[l_j^m(t),l_{j+1}^f(t)),\ {\rm for}\ j=B_{n-1}‚Äô\ldots,J-1\\
\Omega_{(i+1)(j+1)}^c(t)&=&\Omega_{i+1}^m(t)\times \Omega_{j+1}^f(t),\ {\rm for}\ i,j=B_{n-1}‚Äô\ldots,J-1
\end{eqnarray*}
and boundary cohorts:
\begin{eqnarray*}
\Omega_{B_n}^m(t)&=&[0,l_{B_n}^m(t))\\
\Omega_{B_n}^f(t)&=&[0,l_{B_n}^f(t))\\
\Omega_{B_nB_n}^c(t)&=&\Omega_{B_n}^m(t)\times \Omega_{B_n}^f(t)\\
\Omega_{iB_n}^c(t)&=&\Omega_{i}^m(t)\times \Omega_{B_n}^f(t),\ {\rm for}\ i=B_{n-1}‚Äô\ldots,J-1\\
\Omega_{B_nj}^c(t)&=&\Omega_{B_n}^m(t)\times \Omega_{j}^f(t),\ {\rm for}\ j=B_{n-1}‚Äô\ldots,J-1.
\end{eqnarray*}\\
This means that at each internalisation moment $t_n$ we are adding new boundary cohorts to account for the influx of new individuals. At the same time boundary cohorts from the previous time interval $[t_{n-1},t_n)$ became internal ones on $[t_{n},t_{n+1})$. 

It is essential to observe that at time step $t_n$ we deal with $J-B_0+n+1$ male and female cohorts and  $(J-B_0+n+1)^2$ couple cohorts.

For $t_n\leq t<t_{n+1}$ masses and their locations are defined as follows:

\begin{equation}\label{m_mass_location}
\left\{
\begin{array}{rcl}
m_i^m(t)&=&\int_{\Omega_i^m(t)}u^m(t,x)dx,\ i=B_n,\ldots,J\\
x_{B_n}^m(t)&=&
\left\{
\begin{array}{l}
0\ {\rm if}\ m_{B_n}^m(t)=0,\\
\frac{\Pi_{B_n}^m(t)}{m_{B_n}^m(t)},\ {\rm otherwise}\\
\end{array}
\right.\\
x_i^m(t)&=&\frac{1}{m_i^m(t)}\int_{\Omega_i^m(t)}xu^m(t,x)dx,\ i=B_n+1,\ldots,J\\
\end{array}
\right.
\end{equation}
\begin{equation}\label{f_mass_location}
\left\{
\begin{array}{rcl}
m_j^f(t)&=&\int_{\Omega_j^f(t)}u^f(t,y)dy,\ j=B_n,\ldots,J\\
y_{B_n}^f(t)&=&
\left\{
\begin{array}{l}
0\ {\rm if}\ m_{B_n}^f(t)=0,\\
\frac{\Pi_{B_n}^f(t)}{m_{B_n}^f(t)},\ {\rm otherwise}\\
\end{array}
\right.\\
y_j^f(t)&=&\frac{1}{m_j^f(t)}\int_{\Omega_j^m(t)}yu^f(t,y)dy,\ j=B_n+1,\ldots,J\\
\end{array}
\right.
\end{equation}
\begin{equation}\label{c_mass_location}
\left\{
\begin{array}{rcl}
m_{ij}^c(t)&=&\int_{\Omega_{ij}^c(t)}u^c(t,x,y)dxdy,\ i,j=B_n,\ldots,J\\
(x_i^c(t),y_j^c(t))&=&\frac{1}{m_{ij}^c(t)}\int_{\Omega_{ij}^c(t)}(x,y)u^c(t,x,y)dxdy, i,j=B_n,\ldots,J,\\
\end{array}
\right.
\end{equation}\\
where $\Pi_{B_n}^m(t)=\int_{\Omega_{B_n}^m(t)}xu^m(t,x)dx$ and $\Pi_{B_n}^f(t)=\int_{\Omega_{B_n}^f(t)}yu^f(t,y)dy$.\\

Let us notice that we deal with 3 unknowns, thus 3 equations, for each male and female cohort and with 2 unknowns for each couple cohort. Knowing the amount of cohorts explicitly we deduce that at each time step $t_n$ we have $3(J-B_0+n+1)$ equations for males and females and $2(J-B_0+n+1)^2$ equations for couples. Consequently, at $t_n\leq t<t_{n+1}$ we have $8+2(J-B_0+n)$ equations more then we had at $t_{n-1}\leq t<t_{n}$. 

As our aim is to track the dynamics of the above functions describing masses and their locations, we wish to derive the appropriate system of ordinary differential equations. To shorten the notation, from now on, we assume that $t\in[t_n,t_{n+1})$ and we write $B$ instead of $B_n$. Let us start with formulating ODEs for masses and locations of male individuals and symmetric results for female individuals and then move to the more demanding case of the couples. We start with differentiating $m_i^m(t)$, $x_i^m(t)$, $i=B+1,\ldots,J$, and boundary functions $m_B^m(t)$ and $\Pi_B^m(t)$:
\begin{eqnarray*}
\frac{d}{dt}m_i^m(t)&=&\int_{\Omega_i^m(t)}\partial_tu^m(t,x)dx+\frac{d}{dt}l_{i}(t)u^m(t,l_{i}(t))-\frac{d}{dt}l_{i-1}(t)u^m(t,l_{i-1}(t))\\
&=&\int_{\Omega_i^m(t)}\partial_tu^m(t,x)dx+u^m(t,l_{i}(t))-u^m(t,l_{i-1}(t))\\
&=&\int_{\Omega_i^m(t)}\partial_tu^m(t,x)dx+\int_{\Omega_i^m(t)}\partial_xu^m(t,x)dx\\
&=&-\int_{\Omega_i^m(t)}c^m(t,x)u^m(t,x)dx,
\end{eqnarray*}

\begin{eqnarray*}
\frac{d}{dt}m_B^m(t)&=&\int_{\Omega_B^m(t)}\partial_tu^m(t,x)dx+\frac{d}{dt}l_{B}(t)u^m(t,l_{B}(t))\\
&=&\int_{\Omega_B^m(t)}\partial_tu^m(t,x)dx+\frac{d}{dt}l_{B}(t)u^m(t,l_{B}(t))-u^m(t,0)+u^m(t,0)\\
&=&\int_{\Omega_B^m(t)}\partial_tu^m(t,x)dx+\int_{\Omega_B^m(t)}\partial_xu^m(t,x)dx+\int_{\bR^2_+}\beta^m(t,x,y)u^c(t,x,y)dxdy,\\
&=&-\int_{\Omega_B^m(t)}c^m(t,x)u^m(t,x)dx+\int_{\bR^2_+}\beta^m(t,x,y)u^c(t,x,y)dxdy.
\end{eqnarray*}

To differentiate the location function $x_i^m(t)$, for $i=B+1,\ldots,J$, we check how the first moment $\bar{x}_i^m(t)=\int_{\Omega_i^m(t)}xu^m(t,x)dx$ evolves.

\begin{eqnarray*}
\frac{d}{dt}\bar{x}_i^m(t)&=&\int_{\Omega_i^m(t)}x\partial_tu^m(t,x)dx+l_{i}(t)u^m(t,l_{i}(t))-l_{i-1}(t)u^m(t,l_{i-1}(t))\\&=&\int_{\Omega_i^m(t)}x\partial_tu^m(t,x)dx+\int_{\Omega_i^m(t)}\partial_x(xu^m(t,x))dx\\
&=&\int_{\Omega_i^m(t)}x\partial_tu^m(t,x)dx+\int_{\Omega_i^m(t)}x\partial_xu^m(t,x)dx+\int_{\Omega_i^m(t)}u^m(t,x)dx\\
&=&-\int_{\Omega_i^m(t)}xc^m(t,x)u^m(t,x)dx+\int_{\Omega_i^m(t)}u^m(t,x)dx.\\
\end{eqnarray*}

Having $\frac{d}{dt}\bar{x}_i^m(t)$ calculated, it is easy to differentiate the location functions for the internal cohorts:

\begin{eqnarray*}
\frac{d}{dt}x_i^m(t)&=&\frac{d}{dt}\left[\frac{\bar{x}_i^m(t)}{m_i^m(t)}\right]=\frac{\frac{d}{dt}\bar{x}_i^m(t)m_i^m(t)}{\left[m_i^m(t)\right]^2}-\frac{\bar{x}_i^m(t)\frac{d}{dt}m_i^m(t)}{\left[m_i^m(t)\right]^2}\\
&=&\frac{-1}{m_i^m(t)}\int_{\Omega_i^m(t)}xc^m(t,x)u^m(t,x)dx+\frac{1}{m_i^m(t)}\int_{\Omega_i^m(t)}u^m(t,x)dx\\
&&+\frac{1}{m_i^m(t)}\underbrace{\frac{1}{m_i^m(t)}\int_{\Omega_i^m(t)}xu(t,x)dx}_{=x_i^m(t)}\int_{\Omega_i^m(t)}c^m(t,x)u^m(t,x)dx\\
&=&\frac{1}{m_i^m(t)}\int_{\Omega_i^m(t)}[x_i^m(t)-x]c^m(t,x)u^m(t,x)dx+\frac{1}{m_i^m(t)}\int_{\Omega_i^m(t)}u^m(t,x)dx.\\
\end{eqnarray*}

\noindent We also need to differentiate the first moment, $\Pi_B^m(t)$, in the boundary cohort:
\begin{eqnarray*}
\frac{d}{dt}\Pi_B^m(t)&=&\frac{d}{dt}\int_{\Omega_B(t)}xu^m(t,x)dx\\
&=&\int_{\Omega_B^m(t)}x\partial_tu^m(t,x)dx+l_{B}(t)u^m(t,l_{B}(t))-0\cdot u^m(t,0)+0\cdot u^m(t,0)\\
&=&\int_{\Omega_B^m(t)}x\partial_tu^m(t,x)dx+\int_{\Omega_B^m(t)}\partial_x(xu^m(t,x))dx,\\
&=&\int_{\Omega_B^m(t)}x\partial_tu^m(t,x)dx+\int_{\Omega_B^m(t)}x\partial_xu^m(t,x)dx+\int_{\Omega_B^m(t)}u^m(t,x)dx,\\
&=&-\int_{\Omega_B^m(t)}xc^m(t,x)u^m(t,x)dx+\int_{\Omega_B^m(t)}u^m(t,x)dx.\\
\end{eqnarray*}
Because we need to obtain the closed form of the schemes, we use the following  facts that hold for $i=B+1,\ldots,J$:

\begin{eqnarray}
\int_{\Omega_i^m(t)}(x_i^m(t)-x)u^m(t,x)dx&=&0,\label{closed_form_1}\\
\int_{\Omega_i^m(t)}f(x)u^m(t,x)dx&=&f(x_i^m(t))m_i^m(t)+\mathcal{O}\left(\max_{x\in \Omega_i^m(t)} |x-x_i^m(t)|^2\right),\label{closed_form_2}\\ 
&&{\rm for}\ f\in\bf{C^2}(\bR_+).\nonumber
\end{eqnarray}

While the first one can be easily justified based on the definition of $x_i^m(t)$
$$
\begin{array}{rcl}\vspace{2mm}
\int_{\Omega_i^m(t)}[x_i^m(t)-x]u^m(t,x)dx&=&x_i^m(t)\int_{\Omega_i(t)}u^m(t,x)dx-\int_{\Omega_i(t)}xu^m(t,x)dx\\
 &= &x_i^m(t)m_i^m(t)-m_i^m(t)x_i^m(t)=0,\\
\end{array}
$$
the proof of the second one requires equation (\ref{closed_form_1}) and the first order Taylor approximation
\begin{eqnarray*}
\int_{\Omega_i^m(t)}f(x)u^m(t,x)dx&=&\int_{\Omega_i^m(t)}f(x_i^m(t))u^m(t,x)dx\\
  &&+\frac{d}{dx}f(x_i^m(t))\int_{\Omega_i^m(t)}[x-x_i^m(t)]u^m(t,x)dx\\
  &&+\int_{\Omega_i^m(t)}\mathcal{O}\left(|x-x_i^m(t)|^2\right)u^m(t,x)dx\\
  &=&f(x_i^m(t))m_i^m(t)+\|u^m(t,\cdot)\|_{\mathbf{L^1}(\Omega_i^m(t))}\cdot\mathcal{O}\left(\max_{x\in \Omega_i^m(t)}|x-x_i^m(t)|^2\right).
\end{eqnarray*}\\

Furthermore, we observe that
\begin{eqnarray}\label{closed_form_3}
\int_{\Omega_{ij}^c(t)}
\left(
\begin{array}{c}\vspace{1.5mm}
x-x_{ij}^{c}(t)\\
y-y_{ij}^c(t)\\
\end{array}
\right)
u^c(t,x,y)dxdy&=&
\left(
\begin{array}{c}\vspace{1.5mm}
0\\
0\\
\end{array}
\right),\\ \label{closed_form_4}
\int_{\Omega_{ij}^c(t)}
\left(
\begin{array}{c}\vspace{1.5mm}
f_1(x,y)\\
f_2(x,y)\\
\end{array}
\right)u^c(t,x,y)dxdy&=&
\left(
\begin{array}{c}\vspace{1.5mm}
f_1(x_{ij}^c(t),y_{ij}^c(t))\\
f_2(x_{ij}^c(t),y_{ij}^c(t))\\
\end{array}
\right)m_{ij}^c(t) \\ 
&&+\mathcal{O}\left(\max_{(x,y)\in \Omega_{ij}^c(t)}|(x_{ij}^c(t)-x,y_{ij}^c(t)-y)|^2\right),\ {\rm for}\ f_1,f_2\in\bf{C^2}(\bR_+^2).\nonumber
\end{eqnarray}

Equality  (\ref{closed_form_4}) holds because

\begin{eqnarray*}
\int_{\Omega_{ij}^c(t)}
\left(
\begin{array}{c}\vspace{1.5mm}
f_1(x,y)\\
f_2(x,y)\\
\end{array}
\right)
u^c(t,x,y)dxdy&=&\int_{\Omega_{ij}^c(t)}
\left(
\begin{array}{c}\vspace{1.5mm}
f_1(x_{ij}^c(t),y_{ij}^c(t))\\
f_2(x_{ij}^c(t),y_{ij}^c(t))\\
\end{array}
\right)
u^c(t,x,y)dxdy\\
&+&\int_{\Omega_{ij}^c(t)}
\left[
\begin{array}{cc}\vspace{1.5mm}
\frac{\partial f_1}{\partial x} & \frac{\partial f_1}{\partial y}\\ 
\frac{\partial f_2}{\partial x} & \frac{\partial f_2}{\partial y}\\
\end{array}
\right]
(x_{ij}^c(t),y_{ij}^c(t))
\left(
\begin{array}{c}\vspace{1.5mm}
x-x_{ij}^c(t)\\
y-y_{ij}^c(t)\\
\end{array}
\right)
u^c(t,x,y)dxdy\\
&+&\int_{\Omega_{ij}^c(t)} \mathcal{O}\left(\max_{(x,y)\in \Omega_{ij}^c(t) }|(x_{ij}^c(t)-x,y_{ij}^c(t)-y)|^2\right)u^c(t,x,y)dxdy\\
&=&
\left(
\begin{array}{c}\vspace{1.5mm}
f_1(x_{ij}^c(t),y_{ij}^c(t))\\
f_2(x_{ij}^c(t),y_{ij}^c(t))\\
\end{array}
\right)
m_{ij}^c(t)\\
&+&\mathcal{O}\left(\max_{(x,y)\in \Omega_{ij}^c(t) }|(x_{ij}^c(t)-x,y_{ij}^c(t)-y)|^2\right)\cdot\|u^c(t,\cdot,\cdot)\|_{\mathbf{L_1}(\Omega_{ij}^c(t))}\\
\end{eqnarray*}

Now we are ready to derive the system of equations for the male population:

\begin{eqnarray}\nonumber
\frac{d}{dt}m_i^m(t)&=&-\int_{\Omega_i^m(t)}c^m(t,x)u^m(t,x)dx\\ \label{closed_form_m_1}
&=&-c^m(t,x_i^m(t))m_i^m(t)+\mathcal{O}\left(\max_{x\in \Omega_i^m(t)} |x-x_i^m(t)|^2\right),\\\nonumber
\end{eqnarray}
\begin{eqnarray}\nonumber
\frac{d}{dt}x_i^m(t)&=&\frac{1}{m_i^m(t)}\int_{\Omega_i^m(t)}[x_i^m(t)-x]c^m(t,x)u^m(t,x)dx\\ \label{closed_form_m_2}
&&+\frac{1}{m_i^m(t)}\int_{\Omega_i^m(t)}u^m(t,x)dx\\ \nonumber
&=&\frac{1}{m_i^m(t)}[x_i^m(t)-x_i^m(t)]c^m(t,x_i^m(t))m_i^m(t)\\ \nonumber
&&+\frac{1}{m_i^m(t)}m_i^m(t)+\mathcal{O}\left(\max_{x\in \Omega_i^m(t)} |x-x_i^m(t)|^2\right),\\ \nonumber
&=&1+\mathcal{O}\left(\max_{x\in \Omega_i^m(t)}|x-x_i^m(t)|^2\right),\\ \nonumber
\end{eqnarray}

\begin{eqnarray}\nonumber
\frac{d}{dt}m_B^m(t)&=&-\int_{\Omega_B^m(t)}c^m(t,x)u^m(t,x)dx+\int_{\bR^2_+}\beta^m(t,x,y)u^c(t,x,y)dxdy\\ \nonumber
&=&-\int_{\Omega_B^m(t)}c^m(t,0)u^m(t,x)dx-\int_{\Omega_B^m(t)}\partial_xc^m(t,0)(x-0)u^m(t,x)dx\\ \nonumber
&&+\int_{\Omega_B^m(t)}\mathcal{O}\left(|x|^2\right)u^m(t,x)dx+\int_{\bR^2_+}\beta^m(t,x,y)u^c(t,x,y)dxdy\\ \nonumber
&=&-c^m(t,0)m_B^m(t)-\partial_xc^m(t,0)\Pi_B^m(t)+ \mathcal{O}\left(\max_{x\in \Omega_B^m(t)} |x|^2\right)\\ \label{closed_form_m_3}
&&+\sum_{i,j=B}^J\beta^m(t,x_{ij}^c(t),y_{ij}^c(t))m_{ij}^c(t)\\ \nonumber
&&+\sum_{i,j=B}^J \|u^c(t,\cdot,\cdot)\|_{\mathbf{L^1}(\Omega_{ij}^c(t))}\mathcal{O}\left(\max_{(x,y)\in \Omega_{ij}^c(t) }|(x_{ij}^c(t)-x,y_{ij}^c(t)-y)|^2\right)\\ \nonumber
&=&-c^m(t,0)m_B^m(t)-\partial_xc^m(t,0)\Pi_B^m(t)+ \mathcal{O}\left(\max_{x\in \Omega_B^m(t)} |x|^2\right)\\ \nonumber
&&+\sum_{i,j=B}^J\beta^m(t,x_{ij}^c(t),y_{ij}^c(t))m_{ij}^c(t)\\ \nonumber
&&+\max_{i,j=B,\ldots,J}\mathcal{O}\left(\max_{(x,y)\in \Omega_{ij}^c(t) }|(x_{ij}^c(t)-x,y_{ij}^c(t)-y)|^2\right)\nonumber
\end{eqnarray}
and finally
\begin{eqnarray}\nonumber
\frac{d}{dt}\Pi_B^m(t)&=&-\int_{\Omega_B^m(t)}xc^m(t,x)u^m(t,x)dx+\int_{\Omega_B^m(t)}u^m(t,x)dx,\\ \label{closed_form_m_4} 
&=&m_B^m(t)-c^m(t,0)\Pi_B^m(t)+\mathcal{O}\left(\max_{x\in \Omega_B^m(t)}|x|^2\right).\\  \nonumber
\end{eqnarray}
As we mentioned before, one can easily conclude that the system of ODEs for the female population is of similar form, namely
\begin{eqnarray}\nonumber
\frac{d}{dt}m_j^f(t)&=&-\int_{\Omega_j^f(t)}c^f(t,y)u^f(t,y)dy\\ \label{closed_form_f_1}
&=&-c^f(t,y_j^f(t))m_j^f(t)+\mathcal{O}\left(\max_{y\in \Omega_j^f(t)}|y-y_j^f(t)|^2\right),\\\nonumber
\end{eqnarray}
\begin{eqnarray}\nonumber
\frac{d}{dt}y_j^f(t)&=&\frac{1}{m_j^f(t)}\int_{\Omega_j^f(t)}[y_j^f(t)-y]c^f(t,y)u^f(t,y)dy\\ \label{closed_form_f_2}
&&+\frac{1}{m_j^f(t)}\int_{\Omega_j^f(t)}u^f(t,y)dy\\ \nonumber
&=&\frac{1}{m_j^f(t)}[y_j^f(t)-y_j^f(t)]c^f(t,y_j^f(t))m_j^f(t)\\ \nonumber
&&+\frac{1}{m_j^f(t)}m_j^f(t)+\mathcal{O}\left(\max_{y\in \Omega_j^f(t)}|y-y_j^f(t)|^2\right),\\ \nonumber
&=&1+\mathcal{O}\left(\max_{y\in \Omega_j^f(t)}|y-y_j^f(t)|^2\right),\\ \nonumber
\end{eqnarray}

\begin{eqnarray}\nonumber
\frac{d}{dt}m_B^f(t)&=&-\int_{\Omega_B^f(t)}c^f(t,y)u^f(t,y)dy+\int_{\bR^2_+}\beta^f(t,x,y)u^c(t,x,y)dxdy\\ \nonumber
&=&-\int_{\Omega_B^f(t)}c^f(t,0)u^f(t,y)dy-\int_{\Omega_B^f(t)}\partial_xc^f(t,0)(y-0)u^f(t,y)dy\\ \nonumber
&&+\int_{\Omega_B^f(t)}\mathcal{O}\left(|y|^2\right)u^f(t,y)dy+\int_{\bR^2_+}\beta^f(t,x,y)u^c(t,x,y)dxdy\\ \nonumber
&=&-c^f(t,0)m_B^f(t)-\partial_xc^f(t,0)\Pi_B^f(t)+\mathcal{O}\left(\max_{y\in \Omega_B^f(t)}|y|^2\right)\\ \label{closed_form_f_3}
&&+\sum_{i,j=B}^J\beta^f(t,x_{ij}^c(t),y_{ij}^c(t))m_{ij}^c(t)\\ \nonumber
&&+\sum_{i,j=B}^J \|u^c(t,\cdot,\cdot)\|_{\mathbf{L^1}(\Omega_{ij}^c(t))}\mathcal{O}\left(\max_{(x,y)\in \Omega_{ij}^c(t) }|(x_{ij}^c(t)-x,y_{ij}^c(t)-y)|^2\right)\\ \nonumber
&=&-c^f(t,0)m_B^f(t)-\partial_xc^f(t,0)\Pi_B^f(t)+ \mathcal{O}\left(\max_{y\in \Omega_B^f(t)} |y|^2\right)\\ \nonumber
&&+\sum_{i,j=B}^J\beta^f(t,x_{ij}^c(t),y_{ij}^c(t))m_{ij}^c(t)\\ \nonumber
&&+\max_{i,j=B,\ldots,J}\mathcal{O}\left(\max_{(x,y)\in \Omega_{ij}^c(t) }|(x_{ij}^c(t)-x,y_{ij}^c(t)-y)|^2\right)\nonumber
\end{eqnarray}
and finally
\begin{eqnarray}\nonumber
\frac{d}{dt}\Pi_B^f(t)&=&-\int_{\Omega_B^f(t)}yc^f(t,y)u^f(t,y)dy+\int_{\Omega_B^f(t)}u^f(t,y)dy,\\ \label{closed_form_f_4} 
&=&m_B^f(t)-c^f(t,0)\Pi_B^f(t)+\mathcal{O}\left(\max_{y\in \Omega_B^f(t)}|y|^2\right).\\  \nonumber
\end{eqnarray}
To obtain the system of equations for dynamics of populations of couples we need to apply a similar treatment for masses and their locations. 
\begin{eqnarray*}
\frac{d}{dt}m_{ij}^c(t)&=&\frac{d}{dt}\int_{l_{i-1}^m(t)}^{l_{i}^m(t)}\int_{l_{j-1}^f(t)}^{l_{j}^f(t)}u^c(t,x,y)dydx\\
&=&\int_{l_{i-1}^m(t)}^{l_{i}^m(t)}\frac{d}{dt}\left[ \int_{l_{j-1}^f(t)}^{l_{j}^f(t)}u^c(t,x,y)dy\right]dx\\
&&+\frac{d}{dt}{l_{i}^m(t)}\left[ \int_{l_{j-1}^f(t)}^{l_{j}^f(t)}u^c(t,l_{i}^m(t),y)dy\right]\\
&&-\frac{d}{dt}{l_{i-1}^m(t)}\left[ \int_{l_{j-1}^f(t)}^{l_{j}^f(t)}u^c(t,l_{i-1}^m(t),y)dy\right]\\
&=&\int_{l_{i-1}^m(t)}^{l_{i}^m(t)}\Bigg\{ \left[ \int_{l_{j-1}^f(t)}^{l_{j}^f(t)}\partial_tu^c(t,x,y)dy\right]\\
&&+\frac{d}{dt}l_{j}^f(t)u^c(t,x,l_{j}^f(t))-\frac{d}{dt}l_{j-1}^f(t)u^c(t,x,l_{j-1}^f(t))\Bigg\}dx\\
&+&\int_{l_{j-1}^f(t)}^{l_{j}^f(t)}u^c(t,l_{i}^m(t),y)dy- \int_{l_{j-1}^f(t)}^{l_{j}^f(t)}u^c(t,l_{i-1}^m(t),y)dy\\
&=&\int_{\Omega_{ij}^c(t)}\partial_tu^c(t,x,y)dx+\int_{\Omega_i^m(t)}u^c(t,x,l_{j}^f(t))dx-\int_{\Omega_i^m(t)}u^c(t,x,l_{j-1}^f(t))dx\\
&+&\int_{\Omega_j^f(t)}u^c(t,l_{i}^m(t),y)dy-\int_{\Omega_j^f(t)}u^c(t,l_{i-1}^m(t),y)dy\\
&=&\int_{\Omega_{ij}^c(t)}\partial_tu^c(t,x,y)dx+\int_{\Omega_{ij}^c(t)}\partial_xu^c(t,x,y)dy+\int_{\Omega_{ij}^c(t)}\partial_yu^c(t,x,y)dx\\
&=&\int_{\Omega_{ij}^c(t)}\left[\mathcal{T}(t,x,y)-c^c(t,x,y)u^c(t,x,y)\right]dxdy.
\end{eqnarray*}

As previously, to differentiate the locations function $(x_{ij}^c,y_{ij}^c)(t)=\frac{1}{m_{ij}^c(t)}\int_{\Omega_{ij}(t)}(x,y)u^c(t,x,y)dxdy$ we start with its first moment that is $(\bar{x}_{ij}^c,\bar{y}_{ij}^c)(t)=\int_{\Omega_{ij}(t)}(x,y)u^c(t,x,y)dxdy$

\begin{eqnarray*}
\frac{d}{dt}(\bar{x}_{ij}^c,\bar{y}_{ij}^c)(t)&=&\frac{d}{dt}\int_{l_{i-1}^m(t)}^{l_{i}^m(t)}\int_{l_{j-1}^f(t)}^{l_{j}^f(t)}(x,y)u^c(t,x,y)dydx\\
&=&\int_{l_{i-1}^m(t)}^{l_{i}^m(t)}\frac{d}{dt}\left[ \int_{l_{j-1}^f(t)}^{l_{j}^f(t)}(x,y)u^c(t,x,y)dy\right]dx\\
&+&\frac{d}{dt}{l_{i}^m(t)}\left[ \int_{l_{j-1}^f(t)}^{l_{j}^f(t)}(l_{i}^m(t),y)u^c(t,l_{i}^m(t),y)dy\right]\\
&&-\frac{d}{dt}{l_{i-1}^m(t)}\left[ \int_{l_{j-1}^f(t)}^{l_{j}^f(t)}(l_{i-1}^m(t),y)u^c(t,l_{i-1}^m(t),y)dy\right]
\end{eqnarray*}
\begin{eqnarray*}
&=&\int_{\Omega_{ij}^c(t)}(x,y)\partial_tu^c(t,x,y)dxdy+\int_{\Omega_i^m(t)}(x,l_{j}^f(t))u^c(t,x,l_{j}^f(t))dx\\
&&-\int_{\Omega_i^m(t)}(x,l_{j-1}^f(t))u^c(t,x,l_{j-1}^f(t))dx\\
&+&\int_{\Omega_j^f(t)}(l_{i}^m(t),y)u^c(t,l_{i}^m(t),y)dy-\int_{\Omega_j^f(t)}(l_{i-1}^m(t),y)u^c(t,l_{i-1}^m(t),y)dy\\
&=&\int_{\Omega_{ij}^c(t)}(x,y)\partial_tu^c(t,x,y)dxdy\\
&+&\int_{\Omega_{ij}^c(t)}(1,0)u^c(t,x,y)dxdy+\int_{\Omega_{ij}^c(t)}(x,y)\partial_xu^c(t,x,y)dxdy\\
&+&\int_{\Omega_{ij}^c(t)}(0,1)u^c(t,x,y)dxdy+\int_{\Omega_{ij}^c(t)}(x,y)\partial_yu^c(t,x,y)dxdy\\
&=&\int_{\Omega_{ij}^c(t)}(x,y)\left[\mathcal{T}(t,x,y)-c^c(t,x,y)u^c(t,x,y)\right]dxdy+\int_{\Omega_{ij}^c(t)}(1,1)u^c(t,x,y)dxdy.
\end{eqnarray*}

As the derivation of the closed form for couples goes in line with reasoning for the male and female populations, we apply equations (\ref{closed_form_3}) and  (\ref{closed_form_4}). We omit details concerning approximation of particular functions, because similar calculations were performed in the case of male population.

For a generic marriage function $\mathcal{T}(t,x,y)$ we are not able to obtain the closed form, and most we can conclude is:
\begin{eqnarray*}
\frac{d}{dt}m_{ij}^c(t)&=&-c^c(t,x_{ij}^c(t),y_{ij}^c(t))m_{ij}^c(t)+\int_{\Omega_{ij}^c(t)}\mathcal{T}(t,x,y)dxdy\\
\frac{d}{dt}(\bar{x}_{ij}^c(t),\bar{y}_{ij}^c(t))&=&\int_{\Omega_{ij}^c(t)}(x,y)\mathcal{T}(t,x,y)dxdy\\
&-&(x_{ij}^c(t),y_{ij}^c(y))c^c(t,x_{ij}^c(t),y_{ij}^c(t))m_{ij}^c(t)+(1,1)m_{ij}^c(t)\\
\end{eqnarray*}

Hence, in further investigations we consider a specific marriage function given in equation (\ref{marriage_function}). To obtain the desired form of ODEs describing dynamics of couples, we need to approximate the integral of the underlying function over cohort $\Omega_{ij}^c(t)$:

\begin{equation}\label{spec_marriage_funct}
\int_{\Omega^c_{ij}(t)}F(t,u^m(t,x),u^f(t,x),u^c(t,x,y))dxdy=
\end{equation}
$$
\frac{\int_{\Omega^c_{ij}(t)}\left\{\Theta(x,y)h(x)g(y)\left[u^m(t,x)-\int_0^\infty u^c(t,x,y)dy\right]\left[u^f(t,y)-\int_0^\infty u^c(t,x,y)dx\right]\right\}dxdy}
{\gamma+\int_0^\infty h(x)\left[u^m(t,x)-\int_0^\infty u^c(t,x,y)dy\right]dx+\int_0^\infty g(y)\left[u^f(t,y)-\int_0^\infty u^c(t,x,y)dx\right]dy}
$$

To obtain the approximation of the first order with respect to the structure variables, we apply formulas (\ref{closed_form_1})-(\ref{closed_form_4}). Observe that the denominator of function $F$ does not depend on structure variables and hence, there is no need of integrating it over $\Omega^c_{ij}(t)$. Let us start with the integral of the numerator. It consists of four parts which can be approximated separately:

\begin{eqnarray*}
&\int_{\Omega^c_{ij}(t)}&\Theta(x,y)h(x)g(y)u^m(t,x)u^f(t,y)dxdy=\\
&&\Theta(x_i^m(t),y_j^f(t))h(x_i^m(t))g(y_j^f(t))m_i^m(t)m_j^f(t)\\
&&+\mathcal{O}(\max_{x\in \Omega_i^m(t)}|x-x_i^m(t)|^2)
+\mathcal{O}(\max_{y\in \Omega_j^f(t)}|y-y_{j}^f(t)|^2)\\
&\int_{\Omega^c_{ij}(t)}&\Theta(x,y)h(x)g(y)u^m(t,x)\int_0^\infty u^c(t,z,y)dzdxdy=\\
&&\sum_{v=B}^J\Theta(x_i^m(t),y_{vj}^c(t))h(x_i^m(t))g(y_{vj}^c(t))m_i^m(t)m_{vj}^c(t)+\\
&&+\max_{v=B,\ldots,J}\mathcal{O}\left(\max_{(x,y)\in \Omega_{vj}^c(t) }|(x_{vj}^c(t)-x,y_{vj}^c(t)-y)|^2\right)\\
&\int_{\Omega^c_{ij}(t)}&\Theta(x,y)h(x)g(y)u^f(t,y)\int_0^\infty u^c(t,x,z)dzdxdy=\\
&&\sum_{w=B}^J\Theta(x_{iw}^c(t),y_j^f(t))h(x_{iw}^c(t))g(y_j^f(t))m_j^f(t)m_{iw}^c(t)+\\
&&+\max_{w=B,\ldots,J}\mathcal{O}\left(\max_{(x,y)\in \Omega_{iw}^c(t) }|(x_{iw}^c(t)-x,y_{iw}^c(t)-y)|^2\right)\\
&\int_{\Omega^c_{ij}(t)}&\Theta(x,y)h(x)g(y)\int_0^\infty u^c(t,x,z)dz\int_0^\infty u^c(t,z,y)dzdxdy=\\
&&\sum_{v,w=B}^J\Theta(x_{iw}^c(t),y_{vj}^c(t))h(x_{iw}^c(t))g(y_{vj}^c(t))m_{vj}^c(t)m_{iw}^c(t)\\
&&+\max_{v,w=B,\ldots,J}\mathcal{O}\left(\max_{(x,y)\in \Omega_{vw}^c(t) }|(x_{vw}^c(t)-x,y_{vw}^c(t)-y)|^2\right)\\
\end{eqnarray*}

To give the idea how the approximations are obtained, we present calculations of the third part:
$$
\int_{\Omega^c_{ij}(t)}\Theta(x,y)h(x)g(y)u^f(t,y)\int_0^\infty u^c(t,x,z)dzdxdy=\
$$
$$
\int_{l_{i-1}^m(t)}^{l_i^m(t)}h(x)\int_0^\infty u^c(t,x,z)dz \int_{l_{j-1}^f(t)}^{l_j^f(t)}\Theta(x,y)g(y)u^f(t,y)dydx=
$$
$$
\int_{l_{i-1}^m(t)}^{l_i^m(t)}h(x)\int_0^\infty u^c(t,x,z)dz \left[\Theta(x,y_j^f(t))g(y_j^f(t))m_j^f(t)+\|u^f(t,\cdot)\|_{\mathbf{L^1}(\Omega_j^f(t))}\cdot\mathcal{O}(|y-y_j^f(t)|^2)\right]dx=\\
$$
$$
\sum_{w=B}^J\int_{l_{i-1}^m(t)}^{l_i^m(t)}\int_{l_{w-1}^f(t)}^{l_w^f(t)}\left(h(x)\left[\Theta(x,y_j^f(t))g(y_j^f(t))m_j^f(t)+\|u^f(t,\cdot)\|_{\mathbf{L^1}(\Omega_j^f(t))}\mathcal{O}(|y-y_j^f(t)|^2)\right]\right)u^c(t,x,z)dzdx=\\
$$
$$
\sum_{w=B}^Jh(x_{iw}^c(t))[\Theta(x_{iw}^c(t),y_j^f(t))g(y_j^f(t))m_j^f(t)]m_{iw}^c(t)+\\
$$
$$
\sum_{w=B}^Jh(x_{iw}^c(t))\Theta(x_{iw}^c(t),y_j^f(t))g(y_j^f(t))m_j^f(t)\cdot\|u^c(t,\cdot,\cdot)\|_{\mathbf{L^1}(\Omega_{ij}^c(t))}\mathcal{O}\left(\max_{(x,y)\in \Omega_{iw}^c(t) }|(x_{iw}^c(t)-x,y_{iw}^c(t)-y)|^2\right)+\\
$$
$$
\sum_{w=B}^J\|u^f(t,\cdot)\|_{\mathbf{L^1}(\Omega_j^f(t))}\mathcal{O}(|y-y_j^f(t)|^2)\cdot\|u^c(t,\cdot,\cdot)\|_{\mathbf{L^1}(\Omega_{iw}^c(t))}\mathcal{O}\left(\max_{(x,y)\in \Omega_{iw}^c(t) }|(x_{iw}^c(t)-x,y_{iw}^c(t)-y)|^2\right)=\\
$$
$$
\sum_{w=B}^Jh(x_{iw}^c(t))[\Theta(x_{iw}^c(t),y_j^f(t))g(y_j^f(t))m_j^f(t)]m_{iw}^c(t)+\max_{w=B,\ldots,J}\mathcal{O}\left(\max_{(x,y)\in \Omega_{iw}^c(t) }|(x_{iw}^c(t)-x,y_{iw}^c(t)-y)|^2\right)\\
$$


Let us define
\begin{eqnarray}\nonumber
N_{ij}(t)&=&\Theta(x_i^m(t),y_j^f(t))h(x_i^m(t))g(y_j^f(t))m_i^m(t)m_j^f(t)\\\label{Numerator}
&&+\sum_{v=B}^J\Theta(x_i^m(t),y_{vj}^c(t))h(x_i^m(t))g(y_{vj}^c(t))m_i^m(t)m_{vj}^c(t)\nonumber\\
&&+\sum_{w=B}^J\Theta(x_{iw}^c(t),y_j^f(t))h(x_{iw}^c(t))g(y_j^f(t))m_j^f(t)m_{iw}^c(t)\\ \nonumber\\
&&+\sum_{v,w=B}^J\Theta(x_{iw}^c(t),y_{vj}^c(t))h(x_{iw}^c(t))g(y_{vj}^c(t))m_{vj}^c(t)m_{iw}^c(t)\nonumber
\end{eqnarray}
and observe that $N_{ij}(t)$ denotes the approximation of the numerator  of formula (\ref{spec_marriage_funct}).
We can easily conclude that $\bar{N}_{ij}(t)$, defined below, denotes an approximation the numerator  of formula (\ref{spec_marriage_funct}) multiplied by vector $(x,y)$:
\begin{eqnarray}\label{Numerator_bar}
\bar{N}_{ij}(t)&=&(x_i^m(t),y_j^f(t))\Theta(x_i^m(t),y_j^f(t))h(x_i^m(t))g(y_j^f(t))m_i^m(t)m_j^f(t)\\\nonumber
&&+\sum_{v=B}^J(x_i^m(t),y_{vj}^c(t)) \Theta(x_i^m(t),y_{vj}^c(t))h(x_i^m(t))g(y_{vj}^c(t))m_i^m(t)m_{vj}^c(t)\\\nonumber
&&+\sum_{w=B}^J(x_{iw}^c(t),y_j^f(t))(\Theta(x_{iw}^c(t),y_j^f(t))h(x_{iw}^c(t))g(y_j^f(t))m_j^f(t)m_{iw}^c(t)\\ \nonumber
&&+\sum_{v,w=B}^J(x_{iw}^c(t),y_{vj}^c(t))\Theta(x_{iw}^c(t),y_{vj}^c(t))h(x_{iw}^c(t))g(y_{vj}^c(t))m_{vj}^c(t)m_{iw}^c(t)\nonumber.
\end{eqnarray}
The denominator of formula (\ref{spec_marriage_funct}) is made up of two parts, which we approximate  separately:
$$
\int_0^\infty h(x)\left[u^m(t,x)-\int_0^\infty u^c(t,x,z)dz\right]dx=\sum_{i=B}^Jh(x_i^m(t))m_i^m(t)-\sum_{i,j=B}^Jh(x_{ij}^c(t))m_{ij}^c(t)
$$
$$
+\max_{i=B,\ldots,J}\mathcal{O}\left(\max_{x\in \Omega_{i}^m(t) }|x_{i}^m(t)-x|^2\right)+\max_{i,j=B,\ldots,J}\mathcal{O}\left(\max_{(x,y)\in \Omega_{ij}^c(t) }|(x_{ij}^c(t)-x,y_{ij}^c(t)-y)|^2\right)
$$
$$
\int_0^\infty g(y)\left[u^f(t,y)-\int_0^\infty u^c(t,z,y)dz\right]dy=\sum_{j=B}^Jg(y_j^f(t))m_j^f(t)+\sum_{i,j=B}^Jg(y_{ij}^c(t))m_{ij}^c(t)
$$
$$
+\max_{j=B,\ldots,J}\mathcal{O}\left(\max_{y\in \Omega_{j}^f(t) }|y_{j}^f(t)-y|^2\right)+\max_{i,j=B,\ldots,J}\mathcal{O}\left(\max_{(x,y)\in \Omega_{ij}^c(t) }|(x_{ij}^c(t)-x,y_{ij}^c(t)-y)|^2\right)
$$
Let us write
\begin{equation}\label{Denominator}
D_{ij}(t)=\gamma+\sum_{i=B}^Jh(x_i^m(t))m_i^m(t)-\sum_{i,j=B}^Jh(x_{ij}^c(t))m_{ij}^c(t)+
\end{equation}
$$
+\sum_{i=B}^Jg(y_j^f(t))m_j^f(t)-\sum_{i,j=B}^Jg(y_{ij}^c(t))m_{ij}^c(t),
$$
and understand that $D_{ij}(t)$ is an approximation of the denominator of formula (\ref{spec_marriage_funct}). 
Now $\int_{\Omega^c_{ij}(t)}F(t,u^m(t,x),u^f(t,x),u^c(t,x,y))dxdy$ is approximated by $\frac{N_{ij}(t)}{D_{ij}(t)}$, while  $\int_{\Omega^c_{ij}(t)}(x,y)F(t,u^m(t,x),u^f(t,x),u^c(t,x,y))dxdy$ is approximated by $\frac{\bar{N}_{ij}(t)}{D_{ij}(t)}$.
 
%

Having the approximations, we obtain equations for mass and location of couples:

\begin{eqnarray}\label{couples}
\frac{d}{dt}m_{ij}^c(t)&=&-c^c(t,x_{ij}^c(t),y_{ij}^c(t))m_{ij}^c(t)+\frac{N_{ij}(t)}{D_{ij}(t)}+{\rm h.o.t}\\ \nonumber
\frac{d}{dt}(\bar{x}_{ij}^c(t),\bar{y}_{ij}^c(t))&=&\left[(1,1)-(x_{ij}^c(t),y_{ij}^c(y))c^c(t,x_{ij}^c(t),y_{ij}^c(t))\right]m_{ij}^c(t)\\
&+&\frac{\bar{N}_{ij}(t)}{D_{ij}(t)}+{\rm h.o.t},\\\nonumber
\end{eqnarray}
where h.o.t stands for {\it higher order terms}, which are of second order.\\
After neglecting higher order terms in equations (\ref{closed_form_m_1})-(\ref{closed_form_f_4}), we are able to present the numerical scheme for masses and locations for the males, females and couples:

\begin{equation}\label{EBTeq1}
\left\{
\begin{array}{rcl}\vspace{2mm}
\frac{d}{dt}\hat m_i^m(t)&=&-\hat c^m(t,\hat x_i^m(t))\hat m_i^m(t),\ i=B+1,\ldots,J,\\ \vspace{2mm}
\frac{d}{dt}\hat x_i^m(t)&=&1,\ i=B+1,\ldots,J,\\ \vspace{2mm}
\frac{d}{dt}\hat m_B^m(t)&=&-\hat c^m(t,0)\hat m_B^m(t)-\partial_x\hat c^m(t,0)\hat \Pi_B^m(t)\\
&&+\sum_{i,j=B}^J\beta^m(t,\hat x_{ij}^c(t),\hat y_{ij}^c(t))\hat m_{ij}^c(t),\\ \vspace{2mm}
\frac{d}{dt}\hat\Pi_B^m(t)&=&\hat m_B^m(t)-\hat c^m(t,0)\hat \Pi_B^m(t),\\
\hat x_{B}^m(t)&=&
\left\{
\begin{array}{l}
0\ {\rm if}\ \hat m_{B}^m(t)=0,\\
\frac{\hat \Pi_{B}^m(t)}{\hat m_{B}^m(t)},\ {\rm otherwise},\\
\end{array}
\right.\\
\end{array}
\right.
\end{equation}
\begin{equation}\label{EBTeq2}
\left\{
\begin{array}{rcl}\vspace{2mm}
\frac{d}{dt}\hat m_j^f(t)&=&-\hat c^f(t,\hat y_j^f(t))\hat m_j^f(t),\ j=B+1,\ldots,J,\\ \vspace{2mm}
\frac{d}{dt}\hat y_j^f(t)&=&1,\ j=B+1,\ldots,J,\\ \vspace{2mm}
\frac{d}{dt}\hat m_B^f(t)&=&-\hat c^f(t,0)\hat m_B^f(t)-\partial_y\hat c^f(t,0)\hat \Pi_B^m(t)\\
&&+\sum_{i,j=B}^J\beta^f(t,\hat x_{ij}^c(t),\hat y_{ij}^c(t))\hat m_{ij}^c(t),\\ \vspace{2mm}
\frac{d}{dt}\hat \Pi_B^f(t)&=&\hat m_B^f(t)-\hat c^f(t,0)\hat \Pi_B^f(t),\\
\hat y_{B}^f(t)&=&
\left\{
\begin{array}{l}
0\ {\rm if}\ \hat m_{B}^f(t)=0,\\
\frac{\hat \Pi_{B}^m(t)}{\hat m_{B}^f(t)},\ {\rm otherwise,}\\
\end{array}
\right.\\
\end{array}
\right.
\end{equation}
\begin{equation}\label{EBTeq3}
\left\{
\begin{array}{rcl}\vspace{2mm}
\frac{d}{dt}\hat m_{ij}^c(t)&=&-\hat c^c(t,\hat x_{ij}^c(t),\hat y_{ij}^c(t))\hat m_{ij}^c(t)+\frac{\hat N_{ij}(t)}{\hat D_{ij}(t)}, \ i,j=B,\ldots,J,\\ \vspace{2mm}
\frac{d}{dt}(\hat{\bar{x}}_{ij}^c(t),\hat{\bar{y}}_{ij}^c(t))&=&\left[(1,1)-(\hat x_{ij}^c(t),\hat y_{ij}^c(y))\hat c^c(t,\hat x_{ij}^c(t),\hat y_{ij}^c(t))\right]\hat m_{ij}^c(t)+\frac{\hat{\bar{N}}_{ij}(t)}{\hat D_{ij}(t)},\\
(\hat x_{ij}^c(t),\hat y_{ij}^c(t))&=&
\left\{
\begin{array}{l}
0\ {\rm if}\ \hat m_{ij}^c(t)=0,\\
\frac{(\hat{\bar{x}}_{ij}^c(t),\hat{\bar{y}}_{ij}^c(t))}{\hat m_{ij}^c(t)},\ {\rm otherwise,}\\
\end{array}
\right.\\
\end{array}
\right.
\end{equation}
where $\hat N_{ij}(t)$,  $\hat{\bar{N}}_{ij}(t)$ and $\hat D_{ij}(t)$ are given by (\ref{Numerator}), (\ref{Numerator_bar}) and (\ref{Denominator}), respectively, except that all the locations of the cohorts and their masses are replaced with relevant variables with “hats”. 

Observe that the right hand side of the last but one equation, e.i. for $\frac{d}{dt}(\hat{\bar{x}}_{ij}^c(t),\hat{\bar{y}}_{ij}^c(t))$, consists of terms only of the form $\hat m_{ij}^c(t)G(\hat x_{ij}^c(t))$.\\

Above differential equations need to be supplemented with the initial conditions at the beginning of the interval $[t_n,t_{n+1})$. Let us start with ODEs for the males and females. For equations evolving in the interval $[t_0,t_1)$ the initial conditions for our scheme $\hat x_i^m(0)$, $\hat m_i^m(0)$, $\hat y_j^f(0)$, $\hat m_j^f(0)$, for $i,j=B_0,\ldots,J$, are known based on the initial conditions of the original problem, while $\hat m_{B_0}^m(0)=0$ and $\hat \Pi_{B_0}^m(0)=0$. For the corresponding ODEs evolving in $[t_n,t_{n+1})$, where $0<n<N$, we take $\hat x_i^m(t_n)=\lim_{t\rightarrow t_n^-} \hat x_i^m(t)$, $\hat m_i^m(t_n)=\lim_{t\rightarrow t_n^-} \hat m_i^m(t)$, $\hat y_j^f(t_n)=\lim_{t \rightarrow t_n^-} \hat y_j^f(t)$, $\hat m_i^f(t_n)=\lim_{t\rightarrow t_n^-} \hat m_i^f(t)$, for $i,j=B_n,\ldots,J$ and as previously $\hat m_{B_n}^m(t_n)=\hat \Pi_{B_n}^m(t_n)=0$ and $\hat m_{B_n}^f(t_n)=\hat \Pi_{B_n}^f(t_n)=0$.
 
In the case of couples, equations evolving in the interval $[t_0,t_1)$ are supplemented with the following initial conditions: $\hat m_{B_0B_0}^c(0)=\hat m_{iB_0}^c(0)=\hat m_{B_0j}^c(0)=0$, $(\hat x_{B_0B_0}^c(0),\hat y_{B_0B_0}^c(0))=(0,0)$, 
$(\hat{\bar{x}}_{B_0j}^c(0),\hat{\bar{y}}_{B_0j}^c(0))=(0,0)$, for $j=B_0,\ldots,J$ and $(\hat{\bar{x}}_{iB_0}^c(0),\hat{\bar{y}}_{iB_0}^c(0))=(0,0)$, for $i=B_0,\ldots,J$.

\section{Well-posedness and convergence of the method}
\subsection{Well-posedness of the structured population model}
Existence, uniqueness and Lipschitz continuity with respect to initial data and model parameters of  a nonlinear version of  model \eqref{the_system} were established in ref.~\cite{Ulikowska2012} using the framework of nonnegative Radon measures $({\mathcal M^+}({\bR_{+}})\times {\mathcal M^+}({\bR_{+}})\times{\mathcal M^+}({\bR^2_{+}}),d)$, for $d=d_1+d_1+d_2$ with Lipschitz bounded distance
$$
d_i(\mu_1,\mu_2)=\sup\left\{\int_{\bR^i_{+}}\varphi d(\mu_1-\mu_2):\ \varphi\in\mathbf{C}^1(\bR^i_+;\bR)\ {\rm and}\ \|\varphi\|_{W^{1,\infty}}\leq1\right\},\ i=1,2.
$$
It was shown in ref.~\cite{Ulikowska2012} that, under the assumption that $$c^f,c^m,\beta^f,\beta^m,h,g \in \mathbf{C}([0,T];\mathbf{W}^{1,\infty}(\bR_+))$$ and $$c^c,\Theta \in \mathbf{C}([0,T];\mathbf{W}^{1,\infty}(\bR^2_+)),$$ a Lipschitz semigroup of solutions is generated in the sense of Definition \ref{DefLipschSemiflow}.

\subsection{Solvability of the EBT algorithm}
Proof of solvability of the EBT method,  i.e. existence and uniqueness of solutions of  system \eqref{EBTeq1}-\eqref{EBTeq2},  is not straightforward, because of the specific definition of the dynamics of the cohorts involving quotients of the variables with denominator which is not necessarily separated from zero.
Similar problem for the EBT algorithm for a scalar structured population model was solved in  ref.~\cite{GwJabMarUli14}. The strategy proposed in ref.~\cite{GwJabMarUli14} can be applied also in our case.

It consists in considering the problem in a restricted domain, in which the right-hand side is locally bounded and locally Lipschitz-continuous, which yields local in time existence and uniqueness of solutions. The main difficulty is in showing that the solutions remain in the restricted domain locally in time, i.e.  the values of $(\hat{\bar{x}}_{ij}^c(t),\hat{\bar{y}}_{ij}^c(t))$ remain in a cone
\begin{equation}\label{cone}
(0,0)\leq (\hat{\bar{x}}_{ij}^c(t),\hat{\bar{y}}_{ij}^c(t))\leq (C_1(t),C_2(t))\cdot \hat m_{ij}^c(t),\end{equation}
for certain $C_i(t)$, $i=1,2$,  which is equivalent to proving the following inequalities:
\begin{eqnarray}
0&\leq&\frac{\hat{\bar{x}}_{ij}^c(t)}{\hat m_{ij}^c(t)}, \label{Ineq1}\\
\frac{\hat{\bar{x}}_{ij}^c(t)}{\hat m_{ij}^c(t)}&\leq& C_1(t). \label{Ineq2}
\end{eqnarray}
While the first inequality is obvious, the second one needs justification. Let us remind that $\hat x_{ij}^c(t)=\frac{\hat{\bar{x}}_{ij}^c(t)}{\hat m_{ij}^c(t)}$ and take $C_1(t)=t+C$, where $C$ is a constant. Then, the inequality \eqref{Ineq2} is equivalent to $\hat x_{ij}^c(t)-t\leq C$. 

Equation (\ref{couples}) yields
\begin{eqnarray*}
\frac{d}{dt}\left(\hat x_{ij}^c(t)-t\right)&=&
\frac{1}{\hat m_{ij}^c(t)\hat D_{ij}(t)}\left[\hat{\bar{N}}_{ij}(t)\cdot \left(
\begin{array}{c}\vspace{1.5mm}
1\\
0\\
\end{array}
\right)
-\hat x_{ij}^c(t)\hat N_{ij}(t)\right]
\\
=\frac{1}{\hat m_{ij}^c(t)\hat D_{ij}(t)}&\cdot&
\Bigg[(\hat x_i^m(t)-\hat x_{ij}^c(t))\Theta(\hat x_i^m(t),\hat y_j^f(t))h(\hat x_i^m(t))g(\hat y_j^f(t))\hat m_i^m(t)\hat m_j^f(t)\\
&&+(\hat x_i^m(t)-\hat x_{ij}^c(t))\sum_{v=B}^J \Theta(\hat x_i^m(t),\hat y_{vj}^c(t))h(\hat x_i^m(t))g(\hat y_{vj}^c(t))\hat m_i^m(t)\hat m_{vj}^c(t)\\\nonumber
&&+\sum_{w=B}^J(\hat x_{iw}^c(t)-\hat x_{ij}^c(t))(\Theta(\hat x_{iw}^c(t),\hat y_j^f(t))h(\hat x_{iw}^c(t))g(\hat y_j^f(t))\hat m_j^f(t)\hat m_{iw}^c(t)\\ \nonumber
&&+\sum_{v,w=B}^J(\hat x_{iw}^c(t)-\hat x_{ij}^c(t))\Theta(\hat x_{iw}^c(t),\hat y_{vj}^c(t))h(\hat x_{iw}^c(t))g(\hat y_{vj}^c(t))\hat m_{vj}^c(t)\hat m_{iw}^c(t)\Bigg].
\end{eqnarray*}
Taking $$\hat x_{i\max}^c(t)={\rm argmax}_{w=B,\ldots,J}\hat x_{iw}^c(t),$$ we obtain 

$$\frac{d}{dt}\left(\hat x_{i\max}^c(t)-t\right)\leq 0 \quad \text{for} \quad  \hat x_{i\max}^c(t)\geq \hat x_i^m(t)=\tilde{C}+t.$$

\noindent Note that $\frac{d}{dt}\left(\hat x_{i\max}^c(t)-t\right)=\frac{d}{dt}\left(\hat x_{i\max}^c(t)-\hat x^m_i(t)\right)$ and therefore
$|\hat x_{i\max}^c(t)-\hat x^m_i(t)|^+\le |\hat x_{i\max}^c(0)-\hat x^m_i(0)|^+$. Finaly $0\le \hat x_{ij}^c(t)\le \hat{C}+t$.
For further details of the proof, we refer to ref.~\cite{GwJabMarUli14}.

\subsection{Remarks on the convergence of the method}
The essential step in the proof of the convergence of the numerical method is the stability result for the Lipschitz semigroup presented in ref.~\cite{Ulikowska2012}, where the metric space $(E,\rho)$ from Definition \ref{DefLipschSemiflow} is understood as $({\mathcal M^+}({\bR_{+}})\times {\mathcal M^+}({\bR_{+}})\times{\mathcal M^+}({\bR^2_{+}}),d)$, for $d=d_1+d_1+d_2$ and
$$
d_i(\mu_1,\mu_2)=\sup\left\{\int_{\bR^i_{+}}\varphi d(\mu_1-\mu_2):\ \varphi\in\mathbf{C}^1(\bR^i_+;\bR)\ {\rm and}\ \|\varphi\|_{W^{1,\infty}}\leq1\right\},\ i=1,2.
$$
As reported above,  equation (\ref{the_system}) is associated with a Lipschitz semigroup,  \cite{Ulikowska2012}. Therefore, due to Proposition \ref{tangential_ineq}, the remaining step in the proof of convergence boils down to the estimate of the tangential condition, which is the right hand side of the estimate (\ref{estimate}).  To identify the missing part of the proof, let us define two continuous operators, projection and extension, in the following way:
\begin{eqnarray*}
P&:&({\mathcal M^+}({\bR_{+}})\times{\mathcal M^+}({\bR_{+}})\times{\mathcal M^+}({\bR_{+}^2}),d)\rightarrow \bR^{4(J-B_0+n+1)}\times \bR^{3(J-B_0+n+1)^2}\\
E&:&\bR^{4(J-B_0+n+1)}\times \bR^{3(J-B_0+n+1)^2} \rightarrow({\mathcal M^+}({\bR_{+}})\times{\mathcal M^+}({\bR_{+}})\times{\mathcal M^+}({\bR_{+}^2}),d),\\
\end{eqnarray*}
given by the formula
\begin{eqnarray*}
&&P(u^m(t,\cdot),u^f(t,\cdot),u^c(t,\cdot,\cdot))=
((x_i^m(t),m_i^m(t)),(y_i^f(t),m_i^f(t)),(x_{ij}^c(t),y_{ij}^c(t),m_{ij}^c(t)))\\
&&E((x_i^m(t),m_i^m(t)),(y_i^f(t),m_i^f(t)),(x_{ij}^c(t),y_{ij}^c(t),m_{ij}^c(t)))\\
&&\quad\quad =(\sum_{i=B(t)}^Jm_i^m(t)\delta_{\{x_i^m(t)\}},\sum_{j=B(t)}^Jm_j^f(t)\delta_{\{y_j^f(t)\}},\sum_{i=B(t)}^J\sum_{j=B(t)}^Jm_{ij}^c(t)\delta_{\{x_{ij}^c(t),y_{ij}^c(t)\}}).
\end{eqnarray*}
One observes that $P\circ E={\rm Id}$, while the $E\circ P$ does not need to be the identity.

From the derivation of the method we obtain immediately an estimate for the tangential condition between $P(u^m(t,\cdot),u^f(t,\cdot),u^c(t,\cdot,\cdot))$ and \newline  $(\hat x_i^m(t),\hat m_i^m(t)),(\hat y_i^f(t),\hat m_i^f(t)),(\hat x_{ij}^c(t),\hat y_{ij}^c(t),\hat m_{ij}^c(t))$ in the finite dimensional space of the approximation, with the appropriate accuracy. To present the rigorous proof of the convergence of the method we need the estimate of the condition of tangentiality between $(u^m(t,\cdot),u^f(t,\cdot),u^c(t,\cdot,\cdot))$ and \newline $(\sum_{i=B(t)}^Jm_i^m(t)\delta_{\{x_i^m(t)\}},\sum_{j=B(t)}^Jm_j^f(t)\delta_{\{y_j^f(t)\}},\sum_{i=B(t)}^J\sum_{j=B(t)}^Jm_{ij}^c(t)\delta_{\{x_{ij}^c(t),y_{ij}^c(t)\}})$ in the space of measures equipped with an appropriate distance.
Similar proof was carried out in ref.~\cite{GwJabMarUli14}. Due to the lengthy of  technical calculations needed in this step, a complete proof is deferred to a forthcoming paper.

\section*{Acknowledgments}The project (work of KK) was financed by
The National Center for Science awarded based on the number of decisions
DEC-2012/05/E/ST1/02218. \\
AMC was supported by the Emmy Noether Programme
of the German Research Council (DFG).\\
PG is the coordinator of the International Ph.D. Projects Programme of Foundation for Polish Science operated within the Innovative Economy Operational Programme 2007-2013 (Ph.D. Programme: Mathematical Methods in Natural Sciences).

\bibliographystyle{plain}
\bibliography{GKMC_v1}

\end{document}